\documentclass[12pt]{amsart}
\usepackage{amsmath,amsfonts,mathrsfs,amssymb}
\newtheorem{thm}{Theorem}[section]
\newtheorem{theorem}[thm]{Theorem}
\newtheorem{corollary}[thm]{Corollary}
\newtheorem{lemma}[thm]{Lemma}
\newtheorem{proposition}[thm]{Proposition}

\theoremstyle{definition}
\newtheorem{definition}[thm]{Definition}

\theoremstyle{remark}
\newtheorem{remark}[thm]{Remark}

\numberwithin{equation}{section}

\newcommand\dom{\operatorname{dom}}

\newcommand{\slim}{\operatornamewithlimits{s-lim}}

\newcommand\diag{\operatorname{diag}}
\newcommand\ti{\tilde}

\newcommand\wh{\widehat}
\newcommand\ov{\overline}

\newcommand\bC{{\mathbb C}}
\newcommand\bH{{\mathbb H}}
\newcommand\bN{{\mathbb N}}
\newcommand\bR{{\mathbb R}}
\newcommand\bZ{{\mathbb Z}}

\newcommand\sA{{\mathscr A}}
\newcommand\sH{{\mathscr H}}
\newcommand\sF{{\mathscr F}}
\newcommand\sG{{\mathscr G}}
\newcommand\sI{{\mathscr I}}
\newcommand\sM{{\mathscr M}}
\newcommand\sK{{\mathscr K}}
\newcommand\sL{{\mathscr L}}

\newcommand\fm{{\mathfrak m}}

\newcommand\bc{{\mathbf c}}
\newcommand\bs{{\mathbf s}}

\newcommand\sB{{\mathscr B}}
\newcommand\sP{{\mathscr P}}

\newcommand\bg{{\mathbf g}}
\newcommand\bx{{\mathbf x}}
\newcommand\bu{{\mathbf u}}
\newcommand\bv{{\mathbf v}}
\newcommand\re{{\mathrm e}}

\newcommand\rt{{\mathrm t}}

\newcommand\by{{\mathbf y}}

\newcommand\al{\alpha}
\newcommand\de{\delta}

\newcommand\la{\lambda}
\newcommand\si{\sigma}

\begin{document}

\title[Inverse problems for Dirac operators]%
{Inverse spectral problems for Dirac operators with summable potentials}%
\author[S.~Albeverio \and R.~Hryniv \and Ya.~Mykytyuk]%
{S.~Albeverio \and R.~Hryniv \and Ya.~Mykytyuk}%

\address[S.A.]{Institut f\"ur Angewandte Mathematik, Universit\"at Bonn,
Wegelerstr. 6, \mbox{D--53115}, Bonn, Germany;\ SFB 611, Bonn,
Germany;\ BiBoS, Bielefeld, Germany;\ IZKS; CERFIM, Locarno,
Switzerland;\ \and Accademia di Architettura, Mendrisio,
Switzerland}
\email{albeverio@uni-bonn.de}

\address[R.H. and Ya.M.]{Institute for Applied Problems of Mechanics and Mathematics,
3b~Naukova st., 79601 Lviv, Ukraine \and Lviv National University,
1 Universytetska st., 79602 Lviv, Ukraine}%
\email{rhryniv@iapmm.lviv.ua \and yamykytyuk@yahoo.com}%



\date{\today}%
\dedicatory{Dedicated to B.~M.~Levitan, one of the pioneers of this subject}%
\begin{abstract}
The spectral properties of Dirac operators on $(0,1)$ with
potentials that belong entrywise to $L_p(0,1)$, for some
$p\in[1,\infty)$, are studied. The algorithm of reconstruction of
the potential from two spectra or from one spectrum and the
corresponding norming constants is established, and a complete
solution of the inverse spectral problem is provided.
\end{abstract}

\maketitle

\section{Introduction}\label{sec:intr}

The main aim of the present article is to solve the direct and
inverse spectral problems for one-dimensional Dirac operators on a
finite interval under possibly least restrictive assumptions on
their potentials. Namely, the Dirac operators under consideration
are generated by the differential expressions
\[
    \ell_Q := B\frac{d}{dx} + Q(x)
\]
and some boundary conditions, where
\begin{equation}\label{eq:1.BQ}
    B = \begin{pmatrix} 0 & 1 \\ -1 & 0
        \end{pmatrix},
        \qquad
    Q(x) = \begin{pmatrix} q_1(x) & q_2(x) \\ q_2(x) & -q_1(x)
        \end{pmatrix},
\end{equation}
and $q_1$ and $q_2$ are real-valued functions from $L_p(0,1)$,
$p\in[1,\infty)$. To simplify unessential technicalities, we shall
only consider the boundary conditions that correspond to the
Neumann--Dirichlet and Neumann ones in the case of
Sturm--Liouville equations, although other boundary conditions can
be treated in a similar manner (cf.~the study of Sturm--Liouville
operators with nonsmooth potentials and various boundary
conditions in~\cite{HMtwo,SS}). The corresponding Dirac operators
$\sA_1$ and $\sA_2$ in the Hilbert space~$\bH:=L_2(0,1)\times
L_2(0,1)$ act according to the formula $\sA_j \mathbf u = \ell_Q
\mathbf u$ on the domains
\[
    \dom \sA_j := \{ \mathbf u= (u_1,u_2)^{\mathrm t} \mid
    u_1, u_2 \in \mathrm{AC}(0,1), \ell_Q\mathbf u \in \bH,
    u_2(0)=u_j(1)=0 \}.
\]

It is well known~\cite{LS} that the operators $\sA_1$ and $\sA_2$ are
selfadjoint in~$\bH$ and have simple discrete spectra accumulating at
$-\infty$ and $+\infty$. Our primary goal is two-fold: firstly, to give
a complete description of the spectra of $\sA_1$ and $\sA_2$ for
potentials $Q$ of the form~\eqref{eq:1.BQ} with $q_1,q_2 \in L_p(0,1)$
for some $p\in [1,\infty)$---i.e., to solve the direct spectral
problem,---and, secondly, to give an algorithm of reconstruction of
these operators from their spectra or from one spectrum and the
corresponding norming constants---i.e., to solve the inverse spectral
problem.

Ever since P.~Dirac suggested in 1929 the equation (later named after
him) modelling the evolution of spin-$\tfrac12$ particles in the
relativistic quantum mechanics~\cite{Th}, its range of applicability in
various areas of physics and mathematics has been continuously
expanding. In particular, in 1973 Ablowitz, Kaup, Newell, and
Segur~\cite{AKNS} discovered that the Dirac equation is related to a
nonlinear wave equation (the ``modified Korteweg--de Vries equation", a
member of the AKNS--ZS hierarchy, see~\cite{AKNS2,ZS}) in the same
manner as the Schr\"odinger equation is related to the KdV equations,
and this stimulated the increasing interest in direct and inverse
problems for Dirac operators in both physical and mathematical
literature. Earlier in 1966, Gasymov and Levitan solved the inverse
problems for Dirac operators on $\bR_+$ by using the spectral
function~\cite{GasLsp} and by the scattering phase~\cite{GasLscat}.
Their investigations were continued and further developed in many
directions. The reference list is so vast, that we can only mention
those papers, which, in our opinion, are most pertinent to our topic
and refer the reader to the bibliography cited therein for further
material. The books by Levitan and Sargsjan~\cite{LS} and by
Thaller~\cite{Th} may serve as a good introduction to the (respectively
mathematical and physical part of the) theory of Dirac operators.

The inverse scattering theory was developed for Dirac operators on the
axis in~\cite{Fam,Fro,Gre,HJKS,NPT}, for Dirac systems of order~$2n$ on
semiaxis in~\cite{Gas}, and for more general canonical systems on $\bR$
in~\cite{Sakh1}. The nonselfadjoint case was treated in~\cite{MV} and
nonstationary scattering, including point interactions, in~\cite{Ni}
and~\cite{ANT}, respectively. Reconstruction from the spectral function
on semiaxis was done in~\cite{Sa70} for a general boundary condition at
$x=0$ and in~\cite{She} in the case of an interface condition in an
interior point; the general first order systems in
$L_2(\bR_+,\bC^{2n})$ were recently treated in~\cite{LM}. Inverse
problems in the periodic case were studied in~\cite{KhI,Kor}, and the
Weyl--Titschmarsh $m$-function was used to recover the potential of the
Dirac operator in~\cite{Ar} and of the Dirac systems of order $2n$
in~\cite{CG,GKM,Sakh} (see also the detailed reference lists therein).

The inverse problems for Dirac operators on a finite interval have also
been studied in detail. Reconstruction of a continuous potential from
two spectra was carried out in~\cite{GD}, from one spectrum and the
norming constants (in the presence of a Coulomb-type singularity)
in~\cite{Dz}, and from the spectral function in~\cite{Mam}. Explicit
formulae for solutions (based on the degenerate
Gelfand--Levitan--Marchenko equation) in the case where finitely many
spectral data are perturbed were given in~\cite{DKh}. Uniqueness
results for other types of inverse problems were established---e.g.,
for mixed spectral~\cite{Ho} or interior~\cite{MT} data, nonseparated
boundary conditions~\cite{Nab}, or for the weighted Dirac
equations~\cite{Wa}. Ambartsumyan-type theorems were proved
in~\cite{Ho2} and for the matrix case in~\cite{Ki}. Finally, uniqueness
of the inverse problem for general Dirac-type systems of order~$2n$ was
recently established in~\cite{Mal99a,Mal99}.

We observe that in the above-cited papers the inverse spectral problems
for Dirac operators on a finite interval were considered for continuous
potentials only, which excludes, e.\,g., the important case of
piecewise constant potentials. We remove this restriction by allowing
potentials belonging entrywise to $L_p(0,1)$, $p\in[1,\infty)$, and
completely solve the inverse spectral problem for Dirac operators in
this class (see Theorems~\ref{thm:3.evas}, \ref{thm:3.alas},
\ref{thm:5.inv}, and \ref{thm:5.norm}). The main idea of the proof
rests on the fact that the transformation operators for Dirac operators
under consideration satisfy not only the classical
Gelfand--Levitan--Marchenko equation~\eqref{eq:4.GLM}, but also its
counterpart~\eqref{eq:4.krein}, which was used by Krein~\cite{Krein} in
the study of the inverse problem for impedance Sturm--Liouville
equations (see also~\cite{AHM,Fa}). This ``Krein equation" survives the
passage to the limit in the $L_p$-topology and thus allows us to treat
potentials belonging to $L_p(0,1)$ entrywise.

The paper is organized as follows. In Section~\ref{sec:TO}
transformation operators are constructed and some of their properties
are established. Based on this, in Section~\ref{sec:asymp} we find the
asymptotics of eigenvalues and norming constants for the operators
$\sA_1$ and $\sA_2$. The Gelfand--Levitan--Marchenko and Krein
equations, which relate the spectral data and the transformation
operators, are derived in Section~\ref{sec:GLM}, and the solution of
the inverse spectral problem is given in Section~\ref{sec:inv}.
Finally, two appendices contain some facts related to harmonic analysis
and the factorisation theory in operator algebras.

Throughout the paper, we shall denote by $\langle \cdot, \cdot \rangle$
the scalar product in $\bH$ and by $\sM_2$ the algebra of $2\times2$
matrices with complex entries endowed with the operator norm
$|\,\cdot\,|$ of the Euclidean space~$\bC^2$. Where no confusion
arises, we abbreviate $L_p(0,1)$ to $L_p$ and write $L_p(\sM_2)$ for
the space $L_p\bigl((0,1),\sM_2\bigr)$ of $\sM_2$-valued functions on
$(0,1)$ with complex-valued entries and the norm
\[
    \|V\|_{L_p}:= \Bigl(\int_0^1 |V(t)|^p\,dt\Bigr)^{1/p}.
\]
Also, $(x,y)^\mathrm{t}$ shall stand for the column-vector in
$\bC^2$ with components $x$ and $y$.


\section{Transformation operators}\label{sec:TO}

Assume that $Q\in L_p(\sM_2)$, $p\in[1,\infty)$, is of the
form~\eqref{eq:1.BQ} and denote by $U(\cdot)=U(\cdot,\la)$ the
Cauchy matrix corresponding to the equation $\ell_Q \bu = \la
\bu$. In other words, $U$ is a $2\times 2$ matrix-valued function
satisfying the equation
\begin{equation}\label{eq:2.eq}
    B \frac{dU}{dx} + QU = \la U
\end{equation}
and the initial condition $U(0)=I:=\diag (1,1)$. Denoting by
$\bc(\cdot,\la):=\left(c_1(\cdot,\la),c_2(\cdot,\la)\right)^\mathrm{t}$
and
$\bs(\cdot,\la):=\left(s_1(\cdot,\la),s_2(\cdot,\la)\right)^\mathrm{t}$
the solutions of the equation $\ell_Q \bu = \la \bu$ satisfying the
initial conditions $c_1(0,\la)=s_2(0,\la)=1$ and
$c_2(0,\la)=s_1(0,\la)=0$, we find that
\[
    U(x,\la) = \begin{pmatrix}
        c_1(x,\la)&s_1(x,\la)\\c_2(x,\la)&s_2(x,\la)
        \end{pmatrix}.
\]
Our next aim is to derive an integral representation for $U$ of a
special form.

\begin{theorem}\label{thm:2.Cauchy}
Assume that $Q\in L_p(\sM_2)$, $p\in[1,\infty)$. Then
\begin{equation}\label{eq:2.P}
    U(x,\la) = \re^{-\la xB}
        + \int_0^x \re^{-\la(x-2s)B} P(x,s)\,ds,
\end{equation}
where the matrix-valued function $P=P_Q$ has the following
properties:
\begin{itemize}
 \item [(a)] for every $x\in [0,1]$ the function $P(x,\,\cdot\,)$
    belongs to $L_p(\sM_2)$;
 \item [(b)] the mapping ${\mathcal P}_Q:\,x \mapsto P_Q(x,\,\cdot\,)\in L_p(\sM_2)$ is
 continuous on~$[0,1]$;
 \item [(c)] the function ${\mathcal P}_Q$ depends
 continuously in $\mathrm{C}([0,1],L_p(\sM_2))$ on $Q\in L_p(\sM_2)$.
\end{itemize}
\end{theorem}

\begin{proof}
The standard variation of constant arguments show that $U$
satisfies the equivalent integral equation (recall that $B^2=-I$)
\begin{equation}\label{eq:2.Uinteq}
    U(x,\la) = \re^{-\la x B} + \int_0^x
            \re^{-\la(x-t)B}BQ(t)U(t)\,dt,
\end{equation}
which can be solved by the method of successive approximations.
Namely, with
\begin{equation}\label{eq:2.Un}
    U_0(x):=\re^{-\la xB}  \quad\text{and}\quad
    U_{n+1}(x):= \int_0^x \re^{-\la (x-t)B}\,
            BQ(t) U_n(t)\,dt
    \quad\text{for}\quad n\ge0,
\end{equation}
the solution of~\eqref{eq:2.Uinteq} formally equals
 \(                             
    \sum_{n=0}^\infty U_n.
 \)
Assume that we have proved that
\begin{equation}\label{eq:2.Unsum}
    \sum_{n=0}^\infty \|U_n\|_\infty < \infty,
\end{equation}
where $\|U_n\|_\infty:= \sup_{x\in[0,1]}|U_n(x)|$. Differentiating
then the recurrence relations~\eqref{eq:2.Un}, we find that
\[
    U'_n(x) =  -\la B U_n(x) + BQ(x) U_{n-1}(x),
\]
which in view of~\eqref{eq:2.Unsum} shows that the series
$\sum_{n=0}^\infty U_n$ converges in the topology of the
space~$W^1_p((0,1), \sM_2)$ to some $ \sM_2$-valued function~$V$.
This function~$V$ solves~\eqref{eq:2.eq} and satisfies the initial
condition~$V(0)=I$, and hence it coincides with the Cauchy
matrix~$U$.

To justify~\eqref{eq:2.Unsum}, we use the identity
\begin{equation}\label{eq:2.commut}
    \re^{-\la xB}Q(t)=Q(t)e^{\la x B},
        \qquad x,t\in[0,1],
\end{equation}
in the recurrence relations~\eqref{eq:2.Un} and derive the formula
\begin{equation}\label{eq:3.11}
    U_n(x) = \int_{\Pi_n(x)} \re^{-\la(x-2\xi_n(\rt))B}\,
        BQ(t_1) \cdots BQ(t_n)\, dt_1\dots dt_n,
\end{equation}
in which we have set        
\begin{gather*}
    \Pi_n(x) = \{\rt:= (t_1,\dots,t_n)\in \bR^n \mid
                0 \le t_n \le \dots \le t_1\le x\},\\
        \xi_n(\rt) = \sum_{l=1}^n (-1)^{l+1} t_l.
\end{gather*}
Upon the change of variables
 \(                          
    s = \xi_n(\rt), \ y_l = t_{l+1}, \ l=1,2,\dots,n-1,
 \)
we recast the integral in~\eqref{eq:3.11} as
\[                          
    U_n(x) = \int_0^x \re^{-\la(x-2s)B}\,P_n(x,s)\,ds,
\]
where $P_1(x,s) \equiv B Q(s)$ and, for all $n\in\bN$ and $0\le s
\le x\le 1$,
\begin{equation}\label{eq:2.Pn}
    P_{n+1}(x,s) = \int_{\Pi^*_n(x,s)} BQ(s+\xi_n({\mathrm y}))\,
         BQ(y_1)\,\cdots\,BQ(y_n)\, dy_1 \dots dy_n,
\end{equation}
with
\[
    \Pi^*_n(x,s) =  \{{\mathrm y}=(y_1,\dots,y_n)\in \bR^n \mid
        0 \le y_n \le y_{n-1} \le \dots \le y_1
            \le  s+ \xi_n({\mathrm y}) \le x\}.
\]
For convenience, we extend the functions $P_n$, $n\ge2$, to the whole
square $[0,1]\times[0,1]$ by setting $P_n(x,s)=0$ for $0\le x < s \le
1$.

Using the H\"older inequality and Fubini's theorem, we find that, for
every $n\in\bN$, the function $P_{n+1}(x,\cdot)$ belongs to
$L_p(\sM_2)$ and that
\begin{align*}                          
     \|P_{n+1}(x,\,\cdot\,)\|^p_{L_p} &= \int_0^1 |P_{n+1}(x,s)|^p\,ds\\
        &\le (n!)^{1-p} \int_0^1
            \int_{\Pi^*_{n}(x,s)} |Q(s+\xi_{n-1}({\mathrm y}))|^p\,
                |Q(y_1)|^p \cdots |Q(y_{n})|^p\,
                    dy_1 \dots dy_{n}\,ds\\
        &= (n!)^{1-p} \int_{\Pi_{n+1}(x)}
            |Q(t_1)|^p\cdots |Q(t_{n+1})|^p \,dt_1\dots dt_{n+1}\\
        &= \frac1{(n!)^p (n+1)} \Bigl(\int_0^x |Q|^{p}\Bigr)^{n+1}
            \le \frac{\|Q\|_{L_p}^{(n+1)p}}{(n!)^p}.
\end{align*}
Henceforth with $C:= \max_{x\in[-1,1]}|e^{-\la xB}|$ we
have
\[
    |U_n(x)| \le C \int_0^x |P_n(x,s)|\,ds
        \le C \|P_n(x,\,\cdot\,)\|_{L_p}
        \le C \frac{\|Q\|_{L_p}^n}{(n-1)!},
\]
and~\eqref{eq:2.Unsum} follows.

Moreover, the above inequality implies that the series
 \(                         
    \sum_{n=1}^\infty P_n(x,\,\cdot\,)
 \)
converges in $L_p(\sM_2)$ to some function $P(x,\cdot)$ and yields
 the estimate
\begin{equation}\label{eq:2.Pbound}
    \|P(x,\,\cdot\,)\|_{L_p} \le \sum_{n=1}^\infty \frac{\|Q\|_{L_p}^n}{(n-1)!}
        = \|Q\|_{L_p} \exp\{\|Q\|_{L_p}\}
\end{equation}
for all $x\in [0,1]$. This establishes (a).

Assume that $\tilde Q$ is another potential in $L_p(\sM_2)$ and
denote by $\tilde P_n$ the corresponding functions constructed as
above but for $\tilde Q$ instead of $Q$; then similar calculations
on account of the inequality
\begin{align*}
    \bigl|\prod_{k=1}^n a_k - \prod_{k=1}^n b_k\bigr|^p
        &\le   \Bigl(\sum_{k=1}^n |a_k-b_k|
                \prod_{j\ne k} (|a_j| + |b_j|)\Bigr)^p\\
        &\le n^{p-1} \sum_{k=1}^n |a_k-b_k|^p
                \prod_{j\ne k} (|a_j| + |b_j|)^p
\end{align*}
lead to the estimate
\begin{equation}\label{eq:2.Pncont}
    \|P_{n+1}(x,\cdot) - \tilde P_{n+1}(x,\cdot)\|^p_{L_p}
            \le \Bigl(\frac{n+1}{n!}\Bigr)^p\|Q-\tilde Q\|^p_{L_p}
            \Bigl(\|Q\|_{L_p}+  \|\tilde Q\|_{L_p} \Bigr)^{np}.
\end{equation}
It follows that
\begin{equation}\label{eq:2.Pcont}
 \|P_Q(x,\,\cdot\,) - P_{\tilde Q}(x,\,\cdot\,)\|_{L_p}
    \le (1+2r)\re^{2r}\|Q-\tilde Q\|_{L_p}
\end{equation}
as soon as $r$ is such that $\|Q\|_{L_p},\|\tilde Q\|_{L_p}\le r$.

Observe that if $\tilde Q \in \mathrm{C}([0,1], \sM_2)$, then the
functions $\tilde P_n$, $n\ge 2$, are continuous in the
square~$[0,1]\times[0,1]$, and, moreover,
\[
    \max_{0\le x,s\le1}|\tilde P_{n}(x,s)|
            \le \frac{\|\tilde Q\|_\infty^{n}}{n!},
\]
so that the function
 \(
    [0,1] \ni x \mapsto P_{\tilde Q}(x,\,\cdot\,) \in L_p(\sM_2)
 \)
is continuous. Since the potential $Q\in L_p(\sM_2)$ is the limit
in~$L_p(\sM_2)$ of potentials~$\tilde Q_n \in \mathrm{C}([0,1],
\sM_2)$, estimate~\eqref{eq:2.Pcont} yields both assertions (b) and
(c). The proof is complete.
\end{proof}

\begin{corollary}\label{cor:2.TO}
Assume that $Q\in L_p(\sM_2)$ and set
\[
    P^+ = \sum _{n=1}^\infty P_{2n}, \qquad
    P^- = \sum _{n=1}^\infty P_{2n-1},
\]
where the functions $P_n$ are given by formula~\eqref{eq:2.Pn}.
Set $\bc_0(x,\la):=(\cos\la x, \sin\la x)^\mathrm{t}$ and
\begin{equation}\label{eq:2.RK}
\begin{aligned}
    R(x,t)&=R_Q(x,t) := P^+(x,t) + P^-(x,t) J, \\
    K(x,t)&=K_Q(x,t) := \tfrac12\bigl[
        R(x,\tfrac{x-t}2) + R(x,\tfrac{x+t}2) J\bigr],
\end{aligned}
\end{equation}
where $J=\diag\{1,-1\}$. Then the vector-function
$\bc(\,\cdot\,,\la)$ is given by
\begin{equation}\label{eq:2.c-sol}
    \bc(x,\la) = \bc_0(x,\la) + \int_0^x K(x,t)\bc_0(t,\la)\,dt.
\end{equation}
\end{corollary}

\begin{proof}
Using~\eqref{eq:2.commut} and \eqref{eq:2.Pn}, we conclude that
\begin{align*}
    \re^{-\la(x-2s)B} P_{2n}(x,s) &=
            P_{2n}(x,s)\re^{-\la(x-2s)B},\\
    \re^{-\la(x-2s)B} P_{2n-1}(x,s) &=
            P_{2n-1}(x,s)\re^{\la(x-2s)B}.
\end{align*}
Therefore equality~\eqref{eq:2.P} can be written as
\[
  U(x,\la) = \re^{-\la xB}
        + \int_0^x P^+(x,s)\re^{-\la(x-2s)B} \,ds
        + \int_0^x P^-(x,s)\re^{\la(x-2s)B} \,ds.
\]
Observing that
\[                          
    \re^{-\la xB} = \begin{pmatrix} \cos\la x & -\sin\la x \\
                                    \sin\la x & \cos\la
                                    x\end{pmatrix}
\]
and taking the first column of the above equality, we get
\[
    \bc(x,\la) = \bc_0(x,\la)
        + \int_0^x R(x,s) \bc_0(x-2s,\la)\,ds.
\]
Since
\begin{align*}
    \int_0^{x/2} R(x,s) \bc_0(x-2s,\la)\,ds
        &= \frac12 \int_0^x R(x,\tfrac{x-t}2)\bc_0(t,\la)\,dt\\
    \int_{x/2}^x R(x,s) \bc_0(x-2s,\la)\,ds
        &= \int_{x/2}^x R(x,s) J\bc_0(2s-x,\la)\,ds\\
        &= \frac12 \int_0^x R(x,\tfrac{x+t}2)J\bc_0(t,\la)\,dt,
\end{align*}
the required relation follows.
\end{proof}

Equality~\eqref{eq:2.c-sol} shows that the operator $\sI + \sK$
defined by
\[
    (\sI+\sK)\bu(x) = \bu(x) + \int_0^x K(x,t)\bu(t)\,dt
\]
transforms the solution of the equation $\ell_0 \bu = \la\bu$
(i.e., with a potential $Q$ equal to zero identically) subject to
the initial conditions $u_1(0) = 1$, $u_2(0)=0$ into the solution
of the equation $\ell_Q \bu = \la\bu$ satisfying the same initial
conditions. Denote by $\tilde \sA_Q$ the operator in $\bH$ acting
as $\tilde A_Q \bu = \ell_Q\bu$ on the domain
\[
    \dom(\tilde \sA_Q) :=\{ \bu = (u_1,u_2)^{\mathrm t}\in
        W^1_2(0,1)\times W^1_2(0,1) \mid \ u_2(0)=0 \};
\]
then $\sI + \sK$ is in fact the transformation operator for
$\tilde \sA_Q$ and $\tilde \sA_0$, i.e., $\tilde \sA_Q(\sI + \sK)
= (\sI + \sK)\tilde \sA_0 $, see Theorem~\ref{thm:2.KTO}.

The operator $\sK$ possesses some important properties, which we
now establish. Denote by $G_p(\sM_2)$ the set of measurable
$2\times2$ matrix-valued functions~$K$ on $[0,1]\times[0,1]$
having the property that, for each $x$ and $t$ in $[0,1]$, the
matrix-valued functions $K(x,\cdot)$ and $K(\cdot, t)$ belong to
$L_p(\sM_2)$ and, moreover, the mappings
\[
    [0,1]\ni x \mapsto K(x,\cdot)\in L_p(\sM_2),\qquad
    [0,1]\ni t \mapsto K(\cdot,t)\in L_p(\sM_2)
\]
are continuous (i.e., they coincide a.e. with some continuous
mappings from $[0,1]$ into $L_p(\sM_2)$). The set $G_p(\sM_2)$
becomes a Banach space under the norm
\begin{equation}\label{eq:3.gp-norm}
    \|K\|_{G_p} :=
        \max\bigl\{\max_{x\in[0,1]}\|K(x,\cdot)\|_{L_p},
        \max_{t\in[0,1]}\|K(\cdot,t)\|_{L_p}\bigr\}.
\end{equation}
We also denote by $\sG_p(\sM_2)$ the set of the integral operators
$\sK$ in $\bH$ with kernels $K$ from~$G_p(\sM_2)$. Under the
induced norm
 \(
    \|\sK\|_{\sG_p} := \|K\|_{G_p}
 \),
the set $\sG_p(\sM_2)$ becomes an algebra. The algebra
$\sG_p(\sM_2)$ is continuously embedded into the algebra
$\sB(\bH)$ of all bounded operators in $\bH$ since the functions
$K$ belonging to $G_p(\sM_2)$ have finite Holmgren norm~\cite{HS};
moreover, for $\sK\in\sG_p(\sM_2)$ the inequality
$\|\sK\|_{\sB(\bH)}\le \|\sK\|_{\sG_p}$ holds true.

\begin{theorem}\label{thm:2.KinGp}
Assume that $Q\in L_p(\sM_2)$; then the integral operator
$\sK=\sK_Q$ with kernel $K$ of~\eqref{eq:2.RK} belongs to
$\sG_p(\sM_2)$ and, moreover, the mapping $L_p(\sM_2) \ni Q
\mapsto \sK_Q \in \sG_p(\sM_2)$ is continuous.
\end{theorem}

\begin{proof}
In view of relations~\eqref{eq:2.RK}, we have
\begin{equation}\label{eq:2.K}
    K(x,t) = \frac12 \bigl[
        P^+(x,\tfrac{x-t}2)  + P^-(x,\tfrac{x-t}2)J +
        P^+(x,\tfrac{x+t}2)J + P^-(x,\tfrac{x+t}2) \bigr],
\end{equation}
and hence it suffices to prove the assertions of the theorem for the
operators with kernels $P^+(x,\tfrac{x\pm t}2)$ and $P^-(x,\tfrac{x\pm
t}2)$. Since the proof is analogous for all four functions, we shall
give it for the function~$P^+(x,\tfrac{x-t}2)=:\hat P(x,t)$ only.

It follows from the proof of Theorem~\ref{thm:2.Cauchy} that the
function $P^+$ enjoys the properties (a)--(c) of that theorem. Simple
arguments based on the change of variables justify the validity of the
properties (a)--(c) for the kernel $\hat P$. It thus remains to
establish similar properties of $\hat P$ with respect to the
variable~$t$.

Assume first that $Q$ is continuous. Changing the variables $\eta=
s+\xi({\mathrm y})$, $ \tilde y_1 =y_2$, \dots, $\tilde y_{n-1} =
y_n$ in integral~\eqref{eq:2.Pn}, we arrive at the relation
\[
    P_{n+1}(x,s) = \int_s^x BQ(\eta) P_n(\eta,\eta-s)\,d\eta,
\]
which yields
\[
    P^+(x,s) = \int_s^x BQ(\eta) P^-(\eta,\eta-s)\,d\eta.
\]
It follows that
\begin{align*}
    \|\hat P(\cdot,t)\|^p_{L_p}
        &= \int_t^1 \Bigl|
           \int_{\tfrac{x-t}2}^x BQ(\eta)P^-(\eta,\eta-\tfrac{x-t}2)\,d\eta
                    \Bigr|^p dx\\
        &\le \int_0^1 d\eta |Q(\eta)|^p
            \int_t^{2\eta+t} |P^-(\eta,\eta-\tfrac{x-t}2)|^p dx\\
        &\le \|Q\|^p_{L_p} \max_{\eta\in[0,1]}
            \|P^-(\eta,\cdot)\|^p_{L_p} \le  \|Q\|^{2p}
            \exp\{p\|Q\|_{L_p}\},
\end{align*}
cf.~\eqref{eq:2.Pbound}. If $\tilde Q$ is another continuous
potential, then we find analogously that
\begin{align*}
    \|\hat P_Q(\cdot,t) - \hat P_{\tilde Q}(\cdot,t)\|^p_{L_p}
      &\le 2^{p-1} \|Q-\tilde Q\|^p_{L_p}
            \max_{\eta\in[0,1]}\|P^-_Q(\eta,\cdot)\|^p_{L_p}\\
      &\quad+ 2^{p-1} \|\tilde Q\|^p_{L_p}
            \max_{\eta\in[0,1]}\|P^-_Q(\eta,\cdot)
                -P^-_{\tilde Q}(\eta,\cdot)\|^p_{L_p}.
\end{align*}
Recalling inequality~\eqref{eq:2.Pncont}, we conclude that the function
$t\mapsto \hat P_Q(\cdot, t)$, which belongs to
$\mathrm{C}\bigl([0,1],L_p(\sM_2)\bigr)$ if $Q$ is continuous, depends
therein continuously on $Q\in \mathrm{C}\bigl([0,1],\sM_2\bigr)$ with
respect to the topology of $L_p(\sM_2)$. Since the space
$\mathrm{C}\bigl([0,1],\sM_2\bigr)$ is dense in $L_p(\sM_2)$, we show
by continuity that, for every $Q\in L_p(\sM_2)$, the function $\hat
P_Q(\cdot,t)$ belongs to $L_p(\sM_2)$ for every fixed $t\in[0,1]$, that
the mapping $t\mapsto \hat P_Q(\cdot, t)$ is continuous, and that the
continuous $L_p(\sM_2)$-valued function of $t$,
 \(
    t\mapsto P_Q(\cdot, t)
 \),
depends continuously in $\mathrm{C}\bigl([0,1],L_p(\sM_2)\bigr)$
on $Q\in L_p(\sM_2)$.

Summing up, we have shown that the function $\hat P = \hat P_Q$
belongs to $G_p(\sM_2)$ and depends in $G_p(\sM_2)$ continuously
on $Q\in L_p(\sM_2)$. This establishes the theorem.
\end{proof}

\begin{theorem}\label{thm:2.KTO}
Assume that $Q\in L_p(\sM_2)$ and let $\sK$ be an integral
operator with ker\-nel~$K$ of~\eqref{eq:2.RK}. Then $\sI+\sK$ is
the transformation operator for the pair $\tilde \sA_Q$ and
$\tilde \sA_0$, i.e., $\tilde \sA_Q (\sI + \sK) = (\sI +
\sK)\tilde \sA_0$.
\end{theorem}

\begin{proof}
Since $\sK$ belongs to $\sG_p(\sM_2)$ by Theorem~\ref{thm:2.KinGp}
and its kernel $K$ is lower-diagonal, it follows that $\sK$ is a
Volterra operator in~$\bH$ and hence $\sI+\sK$ is a homeomorphism
of~$\bH$.

Write $\hat \sA_Q := (\sI + \sK)^{-1}\tilde \sA_Q (\sI + \sK)$. In view
of~\eqref{eq:2.c-sol}, $\bc_0(\,\cdot\,,\la)$ is an eigenvector of the
operator $\hat\sA_Q$ corresponding to the eigenvalue $\la$ for every
$\la\in\bC$. Denote by $\sL$ the linear hull of the
system~$\{\bc_0(\,\cdot\,,\la) \mid \la\in\bC\}$; then the restrictions
of the operators $\hat \sA_Q$ and $\tilde \sA_0$ onto $\sL$ coincide.
Since $\sL$ is a core of $\tilde \sA_0$ (see a similar result
in~\cite[Theorem~3.3]{AHM}) and $\hat\sA_Q$ is closed, it follows that
$\tilde \sA_0 \subset \hat\sA_Q$. It remains to observe that
$\hat\sA_Q$ cannot be a proper extension of $\tilde \sA_0$ since
otherwise $\hat\sA_Q$---and $\tilde\sA_Q$ by similarity---would have a
two-dimensional nullspace, which would contradict the uniqueness of
solutions to the equation $\ell_Q\bu = \la\bu$. Thus $\tilde \sA_0 =
\hat\sA_Q$, and the proof is complete.
\end{proof}


\section{Direct spectral problem}\label{sec:asymp}

The aim of this section is to perform the direct spectral analysis
for the Dirac operators~$\sA_1$ and $\sA_2$. The main tool of our
investigations will be the transformation operators constructed in
the previous section.

\begin{theorem}\label{thm:3.evas}
Assume that $Q\in L_p(\sM_2)$; then the eigenvalues $(\la_n)_{n\in
\bZ}$ and $(\mu_n)_{n\in \bZ}$ of $\sA_1$ and $\sA_2$ respectively
can be enumerated so that they satisfy the interlacing condition
\begin{equation}\label{eq:3.interlace}
\la_{n-1} < \mu_n < \la_n, \qquad n\in\bZ,
\end{equation}
and the asymptotics
\begin{equation}\label{eq:3.as-la}
 \begin{aligned}
    \la_n  &= \pi (n+\tfrac12) +  e_n(g_1),\\
    \mu_n  &= \pi n +  e_n(g_2),
 \end{aligned}
\end{equation}
where $g_1,g_2\in L_p$ and $e_n(g):= \int_0^1 \mathrm e^{-2\pi
nix}g(x)\,dx$ are the Fourier coefficients of a function~$g$.
\end{theorem}

\begin{proof}
The fact that the spectra of $\sA_1$ and $\sA_2$ interlace (i.e., that
between two consecutive eigenvalues of one operator, there is exactly
one eigenvalue of the other operator) is well known (see, e.g.,
\cite{GD}). It is easily seen that if we have an enumeration of $\la_n$
and $\mu_n$ obeying~\eqref{eq:3.as-la} for some $g_1,g_2\in L_p$,
then~\eqref{eq:3.interlace} holds for all $n$ with sufficiently
large~$|n|$. Since the two spectra interlace, we can permute a finite
number of indices if necessary in such a way
that~\eqref{eq:3.interlace} becomes valid for all $n\in\bZ$. This
reordering amounts to adding trigonometric polynomials of finite degree
to the functions $g_1$ and $g_2$, and the modified functions $g_1$ and
$g_2$ will remain in $L_p$. Therefore it suffices to
establish~\eqref{eq:3.as-la}.

The equation $\ell_Q \bu = \la \bu$ subject to the initial
conditions $f_1(0)=1$, $f_2(0)=0$ has the solution
\[
    \mathbf c(x,\la) = \mathbf c_0(x,\la)
        + \int_0^x K(x,s) \mathbf c_0(s,\la)\,ds,
\]
where $K=:(k_{jl})_{j,l=1}^2$ is the kernel of the transformation
operator constructed in Corollary~\ref{cor:2.TO}. The numbers $\la_n$
are zeros of the function $c_1(1,\la)$, which, after simple
transformations, takes the form
\begin{equation}\label{eq:3.phiint}
 \begin{aligned}
    c_1(1,\la)
        &= \cos\la + \int_0^1 \bigl[k_{11}(s) \cos(\la s)
                    + k_{12}(s) \sin(\la x)\bigr]\,ds\\
        &= \cos\la + \int_{-1}^1 f_1(s) \mathrm e^{i\la s}\,ds,
 \end{aligned}
\end{equation}
where
\[
    f_1(s) := \begin{cases}
            \tfrac12 \bigl[ k_{11}(1,s) - i k_{12}(1,s)\bigr],
                    &\quad s \ge 0,\\
            \tfrac12 \bigl[ k_{11}(1,-s) +i k_{12}(1,-s)\bigr],
                    &\quad s < 0
            \end{cases}
\]
is a function in~$L_p(-1,1)$.

Analogously the eigenvalues of $\sA_2$ are zeros of the entire
function $c_2(1,\la)$, which has the form
\begin{equation}\label{eq:3.psiint}
    c_2(1,\la) = \sin\la +  \int_{-1}^1 f_2(s) \mathrm e^{i\la s}\,ds
\end{equation}
for some $f_2\in L_p(-1,1)$.

The required asymptotics for zeros of $c_1(1,\la)$ and
$c_2(1,\la)$ follows now from~\cite{HMzero}.
\end{proof}

\begin{definition}\label{def:3.SD}
We denote by $\mathrm{SD}_p$ the set of all pairs
$\{(\la_n)_{n\in\bZ},(\mu_n)_{n\in\bZ}\}$, in which $(\la_n)$ and
$(\mu_n)$ are sequences of real numbers that obey the interlacing
condition~\eqref{eq:3.interlace} and the
asymptotics~\eqref{eq:3.as-la}.
\end{definition}

Theorem~\ref{thm:3.evas} shows that, for any real-valued $Q\in
L_p(\sM_2)$, the spectra of the operators $\sA_1$ and $\sA_2$ form
an element of $\mathrm{SD}_p$. In the reverse direction,
Theorem~\ref{thm:5.inv} claims that any element of $\mathrm{SD}_p$
is composed of the spectra of Dirac operators $\sA_1$ and $\sA_2$
for some $Q\in L_p(\sM_2)$. The reconstruction algorithm uses in
fact the spectrum $(\la_n)$ of $\sA_1$ and the sequence of
corresponding norming constants $\al_n :=
\|\bc(\,\cdot\,,\la_n)\|^{-2}$, whose properties we are going to
study next.

\begin{theorem}\label{thm:3.alas}
Assume that $Q\in L_p(\sM_2)$; then the norming constants
$\al_n=\|\bc(\,\cdot\,,\la_n)\|^{-2}$ have the asymptotics
\begin{equation}\label{eq:3.as-al}
    \al_n = 1 +  e_n(g),
\end{equation}
where $g\in L_p$.
\end{theorem}

Set $\phi(\la):=c_1(1,\la)$ and $\psi(\la):=c_2(1,\la)$. We
observe that the functions $\phi(\la)$ and $\psi(\la)$ can be
reconstructed from their zeros as follows~\cite{HMzero}.

\begin{proposition}\label{pro:3.PhiPsiprod}
The following equalities hold:
\begin{equation}\label{eq:3.PhiPsiprod}
    \phi(\lambda) = \mathrm{V.p.}\hspace{-3pt}
        \prod\limits_{n=-\infty}^{\infty}
        \frac{\la_n -\lambda}{\pi(n+\tfrac12)},     \qquad
    \psi(\lambda) = (\lambda-\mu_0) \mathrm{V.p.}\hspace{-3pt}
        \prod\limits_{n=-\infty}^{\infty}
        \hspace{-7pt}{\vphantom{\prod}}^\prime
        \hspace{5pt}\frac{\mu_n-\lambda}{\pi n}
\end{equation}
(where the prime means that the factor corresponding to the index $n=0$
should be omitted); moreover, the products converge uniformly on
compact sets.
\end{proposition}

It turns out that $\phi$ and $\psi$ (i.e., that the two spectra
$(\la_n)_{n\in\bZ}$ and  $(\mu_n)_{n\in\bZ}$) determine the
norming constants $\al_n$  as follows.

\begin{lemma}\label{lem:3.al}
The norming constants $\al_n$ satisfy the following relation:
\begin{equation}\label{eq:3.al}
    \al_n = - \frac{1}{\dot{\phi}(\la_n)\psi(\la_n)}.
\end{equation}
\end{lemma}

\begin{proof}
Since the system of eigenvectors $(\bc_n)_{n\in {\bZ}}$,
$\bc_n:=\bc(\cdot,\la_n)$, of the operator $\sA_1$ is an
orthogonal basis of $\bH$, we have
\[
    (\sA_1-\la)^{-1}\bu  =  \sum_{n=-\infty}^{\infty}
        \frac{\alpha_n \langle\bu, \bc_n\rangle \bc_n}
        {\la_n -\lambda}.
\]
In particular, the residue of this expression at $\la=\la_n$
equals $-\alpha_n \langle\bu, \bc_n\rangle \bc_n$.

On the other hand, $(\sA_1-\la)^{-1}$ can be calculated as
\begin{align*}
    (\sA_1-\la)^{-1}\bu (x)
        &=  \Bigl[\int_x^1 \langle \bu(t),
            \ov{\bs^+(t,\la)} \rangle_{\bC^2}dt\Bigr]
                \frac{\bc(x,\la)}{W(\la)}\\
        &\quad+ \Bigl[\int_0^x \langle \bu(t),
            \ov{\bc(t,\la)}\rangle_{\bC^2}dt\Bigr]
            \frac{\bs^+(x,\la)}{W(\la)},
\end{align*}
where $\bs^+(\cdot,\la)$ is a solution of the equation $\ell_Q
(\bu) = \la^2\bu$ subject to the terminal conditions $u_1(1)=0$
and $u_2(1)=1$, and $W(\lambda)$ is the Wronskian of the solutions
$\bc(\cdot,\lambda)$ and $\bs^+(\cdot,\lambda)$, i.e.,
\[
 W(\lambda):=   c_1(x,\lambda) s^+_2(x,\lambda)
        -  c_2(x,\lambda) s^+_1(x,\lambda) .
\]
Since the right-hand side of the above equality does not depend on~$x$,
we see that
\[
    W(\lambda)= c_1(1,\lambda) = \phi(\la).
\]
Thus the residue of $(\sA_1-\la)^{-1}\bu$ at $\la=\la_n$ is equal
to
\begin{equation}\label{eq:3.residW}
    \Bigl[\int_x^1 \langle \bu(t),
            \ov{\bs^+(t,\la_n)} \rangle_{\bC^2}dt\Bigr]
        \frac{\bc(x,\la_n)}{\dot{\phi}(\la_n)}
        + \Bigl[\int_0^x \langle \bu(t),
            \ov{\bc(t,\la_n)}\rangle_{\bC^2}dt\Bigr]
        \frac{\bs^+(x,\la_n)}{\dot{\phi}(\la_n)}.
\end{equation}
We observe now that the vector-functions $\bs^+(x,\la_n)$ and
$\bc(x,\la_n)$ are collinear, namely,
\[
    \bs^+(\cdot,\la_n)
        = \frac{s^+_2(1,\la_n)}{c_2(1,\la_n)}\bc(\cdot,\la_n)
        = \frac1{\psi(\la_n)}\bc(\cdot,\la_n),
\]
and thus expression~\eqref{eq:3.residW} simplifies to
$\langle\bu,\bc_n\rangle\bc_n/(\dot{\phi}(\la_n)\psi(\la_n))$. Equating
this with the above expression for the residue of $(\sA_1-\la)^{-1}\bu$
at $\la=\la_n$, we obtain relation~\eqref{eq:3.al} for $\al_n$. The
proof is complete.
\end{proof}

\begin{proof}[Proof of Theorem~\ref{thm:3.alas}]
In view of equality~\eqref{eq:3.al} and
Propositions~\ref{pro:A.mult} and \ref{pro:A.wiener} it suffices
to prove that the numbers $a_n:=-\dot\phi(\la_n)$ and
$b_n:=\psi(\la_n)$ can be represented as
 $(-1)^{n}(1+\tilde a_n)$ and $(-1)^{n}(1+\tilde b_n)$ respectively,
where $\tilde a_n$ and $\tilde b_n$ are $n$-th Fourier
coef\-fi\-ci\-ents of some functions from $L_p$.

Using formulae~\eqref{eq:3.phiint} and \eqref{eq:3.psiint}, we show
that the numbers $a_n$ and $b_n$ are of the form
\[
     \sin\la_n + \int_{-1}^1 f(s)\re^{i\la_n s}\, ds
\]
with $f(s)=-is f_1(s)\in L_p(-1,1)$ for $a_n$ and $f(s)=f_2(s)\in
L_p(-1,1)$ for $b_n$. Hence it remains to show that
 $(-1)^{n}\sin\la_n = 1 + e_{n}(\tilde g_1)$ for some
$\tilde g_1\in L_p$ and that
\[
    (-1)^{n}\int_{-1}^1 f(s)\re^{i\la_n s}\, ds = e_{n}(\tilde g_2)
\]
for some $\tilde g_2\in L_p$.

Since by Theorem~\ref{thm:3.evas} the eigenvalues $\la_n$ satisfy
the relation $\la_n = \pi (n+\tfrac12) + \ti\la_n$, where
$\ti\la_n=e_{n}(g_1)$ for some function $g_1\in L_p$, we find that
$\sin\la_n = (-1)^n \cos\ti\la_n$, so that
\[
    (-1)^{n} \sin\la_n -1
        = \sum_{k=1}^\infty (-1)^k
        \frac{\ti\la_n^{2k}}{(2k)!}.
\]
Applying Proposition~\ref{pro:A.mult} to the element
$\bx:=(\ti\la_n)_{n\in\bZ}\in X_p$, we see that there exists a
function $\tilde g_1$ in $L_p$ such that the sum on the right-hand
side of the above equality equals~$e_{n}(\tilde g_1)$.

Changing the variables $s\mapsto 1-2t$ in the integral, we get
\begin{align*}
    \int_{-1}^1 f(s)\re^{i\la_n s}\, ds
        &= \int_{0}^1 2f(1-2t)\re^{i\la_n (1-2t)}\, dt\\
        &=  (-1)^{n}\int_{0}^1 \tilde f(t)\re^{i\tilde \la_n (1-2t)}
            \re^{-2 \pi nit}\, dt
\end{align*}
with $\tilde f(t):= 2i f(1-2t)\re^{-\pi it}$. Developing the
function~$\re^{i\ti\la_n (1-2t)}$ into the Taylor series and then
changing summation and integration order (which is allowed in view
of the absolute convergence of the Taylor series and the
integral), we find that
\begin{equation*}
\begin{aligned}
    (-1)^{n}\int_{-1}^1 f(s)\re^{i\la_n s}\, ds
        &  =     \sum_{k=0}^\infty \frac{\tilde \la_n^{k}}{k!}
                    e_{n}(V^{k}\tilde f),
\end{aligned}
\end{equation*}
where $V$ is the operator of multiplication by $i(1-2t)$. In
virtue of Proposition~\ref{pro:A.mult} and the fact that $V$ has
norm~$1$ in $L_p$, the right-hand side of the above equality gives
the $n$-th Fourier coefficient of some function $\tilde g_2$ from
$L_p$, and the proof is complete.
\end{proof}


\section{Derivation of the GLM and the Krein equation}\label{sec:GLM}

Write $\bc_n(\cdot)=\bc(\cdot,\la_n)$ and
$\al_n:=\|\bc_n\|_{\bH}^{-2}$. Since the functions
$\{\bc_n\}_{n\in\bZ}$ form an orthogonal basis of~$\bH$, we have
\[
    \slim_{k\to\infty} \sum_{n=-k}^k \al_n \langle\cdot, \bc_n\rangle\bc_n
    =\sI,
\]
where $\sI$ is the identity operator in $\bH$. On the other hand,
$\bc_n = (\sI+\sK)\bv_n$ with $\bv_n(x):=\bc_0(x,\la_n)=(\cos\la_n
x,\sin\la_n x)^\mathrm{t}$, so that the previous relation can be
rewritten as
\[
    (\sI+\sK) \Bigl[\, \slim_{k\to\infty}
        \sum_{n=-k}^k  \al_n \langle\cdot, \bv_n\rangle \bv_n \Bigr]
        (\sI+\sK)^* = \sI,
\]
which implies
\[
    \slim_{k\to\infty}
        \sum_{n=-k}^k \al_n \langle\cdot, \bv_n\rangle\bv_n
        = (\sI+\sK)^{-1}(\sI+\sK^*)^{-1}.
\]

Set
\begin{equation}\label{eq:4.Hexp}
    H(s):= \mathrm{V.p.}\sum_{n=-\infty}^\infty
        \Bigl(\al_n \re^{-2\la_n sB} - \re^{-\pi (2n+1) sB}\Bigr),
\end{equation}
where for $p=1$ the summation is understood in the Ces\`{a}ro sense,
see Lemma~\ref{lem:5.H}. Observing that
$\langle\,\cdot\,,\bc_{0}(\cdot,\la)\rangle\bc_{0}(\cdot,\la)$ is an
integral operator with kernel
\[
     \tfrac12 \bigl[\re^{-\la (x-t) B}
        + \re^{-\la(x+t)B}J\bigr]
\]
and that
\[
     \slim_{k\to\infty}
        \sum_{n=-k}^k \langle\cdot,\bv_{0,n} \rangle\bv_{0,n} =\sI
\]
with $\bv_{0,n} := \bc_{0}\bigl(\cdot,\pi (n+\tfrac{1}2)\bigr)$,
we conclude that
\begin{equation}\label{eq:4.factor}
    (\sI+\sK)^{-1} (\sI+\sK^*)^{-1} = \sI+\sF,
\end{equation}
where $\sF$ is an integral operator in~$\bH$,
\[
    (\sF \bu)(x) = \int_0^1 F(x,s)\bu(s)\,ds,
\]
with kernel
\begin{equation}\label{eq:4.F}
    F (x,s) := \tfrac12 \bigl[ H(\tfrac{x-t}2) +
                   H(\tfrac{x+t}2)J\bigr].
\end{equation}

Since $\sK^*$ is an integral Volterra operator with upper-diagonal
kernel, the operator $\sI + \sK^*$ is invertible and its inverse can be
written in the form $\sI + \tilde \sK$, where $\tilde \sK$ is an
integral operator in~$\bH$ with upper-diagonal kernel.
By~\eqref{eq:4.factor} one gets
\[
    (\sI+\sK)(\sI+\sF) = (\sI+\sK^*)^{-1}
        = \sI + \tilde \sK;
\]
spelling out this equality in terms of the kernels $K$ and $F$ for
$x>t$, we arrive at the Gelfand--Levitan--Marchenko (GLM) equation
\begin{equation}\label{eq:4.GLM}
    K(x,t) + F(x,t)
        + \int_0^x K(x,s) F(s,t)\,ds = 0.
\end{equation}

If $Q$ is continuous, then such is also $K$, and one has the
formula~\cite[Lemma~12.1.1]{LS}
\begin{equation}\label{eq:4.QKB}
    Q(x) = K(x,x)B - B K(x,x)
\end{equation}
relating the potential and the kernel of the corresponding
transformation operator. This suggests the following algorithm of
solution of the inverse spectral problem: given the spectral data
$\{(\la_n),(\al_n)\}$, one constructs first the kernel $F$
via~\eqref{eq:4.F} and \eqref{eq:4.Hexp}, then solves the GLM
equation~\eqref{eq:4.GLM} for $K$, and, finally, recovers the potential
via~\eqref{eq:4.QKB}. However, if $Q\in L_p(\sM_2)$, then
relation~\eqref{eq:4.QKB} becomes meaningless since, by
Theorem~\ref{thm:2.KinGp}, in this case the kernel $K$ belongs to
$G_p(\sM_2)$ but can be neither continuous nor well defined on subsets
of $[0,1]\times[0,1]$ of Lebesgue measure zero.

It turns out that for $Q\in L_p(\sM_2)$ the restriction $R(x,x)$
of the kernel $R$ of~~\eqref{eq:2.RK} to the diagonal does
determine a matrix-function with entries in $L_p(0,1)$. If $Q$ is
continuous, then~\eqref{eq:4.QKB} together with~\eqref{eq:2.RK}
and the commutator relations
\begin{equation}\label{eq:4.commut}
    BR(x,t) = R(x,t) B, \qquad
    BJ = - JB
\end{equation}
yields the equality
\begin{equation}\label{eq:4.QR}
    Q(x)  =  R(x,x)JB.
\end{equation}
Equation~\eqref{eq:4.QR} retains sense also for $Q\in L_p(\sM_2)$ as an
equality in $L_p(\sM_2)$ and thus can be used to recover~$Q$ in this
situation.

Our next task is to explain how the kernel $R$ can be determined from
the spectral data. We do this by deriving below the Krein
equation~\eqref{eq:4.krein}, an analogue of the GLM equation for
$R$~\cite{Fa,Krein}.

Applying the operator~$B$ to equality~\eqref{eq:4.GLM} from both sides,
we obtain its counterpart,
\begin{equation}\label{eq:4.GLM-B}
    BK(x,t)B + BF(x,t)B
        + \int_0^x BK(x,s) F(s,t)B\,ds = 0.
\end{equation}
With regard to the commutator relations~\eqref{eq:4.commut} and $BH(x)
= H(x)B$, we see that
\begin{align*}
    BK(x,t)B &= \tfrac12 \bigl[
        - R(x,\tfrac{x-t}2) + R(x,\tfrac{x+t}2) J\bigr],\\
    BF(x,t)B &= \tfrac12 \bigl[ -H(\tfrac{x-t}2) +
                   H(\tfrac{x+t}2)\bigr]
\end{align*}
and hence
\begin{align*}
    K(x,t) + BK(x,t)B &= R(x,\tfrac{x+t}2)J,\\
    F(x,t) + BF(x,t)B &= H(\tfrac{x+t}2)J.
\end{align*}
It also follows that
\[
    K(x,s)F(s,t) + BK(x,s)F(s,t)B
        = \tfrac12\bigl[
          R(x,\tfrac{x-s}2) H(\tfrac{s+t}2)J
        + R(x,\tfrac{x+s}2)JH(\tfrac{s-t}2)   \bigr].
\]
Adding now~\eqref{eq:4.GLM} and \eqref{eq:4.GLM-B}, combining the
above formulae in the resulting expression, and using the relation
 \( JH(x) = H(-x)J \),
we arrive at the equation
\[
    R(x,\tfrac{x+t}2) + H(\tfrac{x+t}2) +
        \int_0^x R(x,s) H(\tfrac{x+t}2-s)\,ds =0,
        \qquad 0\le t\le x \le 1.
\]
Subtracting~\eqref{eq:4.GLM-B} from~\eqref{eq:4.GLM} and
performing similar transformations, we arrive at the above formula
with $\tfrac{x+t}2$ replaced by $\tfrac{x-t}2$, and both can now
be combined together to give
\[
    R(x,t) + H(t) + \int_0^x R(x,s) H(t-s)\,ds =0,
        \qquad 0\le t\le x \le 1.
\]
We see that the function $\tilde R(x,t):= R(x,x-t)$ satisfies the
following Krein equation:
\begin{equation}\label{eq:4.krein}
    \tilde R(x,t) + H(x-t) + \int_0^x \tilde R(x,s)H(s-t)\,ds =0.
\end{equation}

We observe that as soon as a kernel $H$ is given by~\eqref{eq:4.Hexp}
with $\la_n$ and $\al_n$ obeying the proper asymptotics (guaranteeing
that the series for $H$ converges in $L_p(\sM_2)$), the integral
operator $\sH$ with kernel $H(x-t)$ belongs to the algebra
$\sG_p(\sM_2)$ introduced in Section~\ref{sec:TO} and $\sI+\sH$ is
positive in $\bH$ (see Lemma~\ref{lem:5.Iholds}). The Krein
equation~\eqref{eq:4.krein} is then uniquely soluble and its solution
belongs to $G_p(\sM_2)$, see Appendix~\ref{app:krein}. In particular,
$\tilde R(x,0)=R(x,x)$ is in $L_p(0,1)$ entrywise indeed.

\begin{remark}\label{rem:4.GLM}
We notice that the GLM equation~\eqref{eq:4.GLM} is the even part
of the Krein equation~\eqref{eq:4.krein} in the sense that if
$\tilde R$ is a solution to~\eqref{eq:4.krein}, then the function
\begin{equation}\label{eq:5.k-and-r}
    K(x,t)
        := \tfrac12 \bigl[\tilde R(x,\tfrac{x-t}{2})
            + \tilde R(x,\tfrac{x+t}{2}J)\bigr]
\end{equation}
solves~\eqref{eq:4.GLM}. Moreover, the condition $\sI+\sH>0$
implies that the operator $\sI + \sF$ is positive in $\bH$ and
thus, in view of the results of Appendix~\ref{app:krein},
guarantees that the GLM equation~\eqref{eq:4.GLM} with $F$
of~\eqref{eq:4.F} is uniquely soluble for $K$ and the solution
belongs to $G_p(\sM_2)$.
\end{remark}


\section{Inverse spectral problem}\label{sec:inv}

The purpose of this section is, firstly, to show by limiting
arguments that formula~\eqref{eq:4.QR} remains valid if the matrix
potential $Q$ belongs to $L_p(\sM_2)$ and, secondly, to justify
the algorithm reconstructing the potential~$Q$ from the spectral
data. Namely, we shall prove the following theorem, which
constitutes the main result of the paper.

\begin{theorem}\label{thm:5.inv}
Assume that $\{(\la_n)_{n\in\bZ},(\mu_n)_{n\in\bZ}\}$ is an
arbitrary element of $\mathrm{SD}_p$, $p\in[1,\infty)$. Then there
exists a unique potential $Q\in L_p(\sM_2)$ such that $(\la_n)$
and $(\mu_n)$ are eigenvalues of the corresponding Dirac operators
$\sA_1$ and $\sA_2$ respectively. The potential $Q$ is equal to
$R(x,x)JB=\tilde R(x,0)JB$, where $\tilde R$ is the solution of
the Krein equation~\eqref{eq:4.krein} with $H$
of~\eqref{eq:4.Hexp} and $\al_n$ given by~\eqref{eq:3.al}.
\end{theorem}

The reconstruction algorithm proceeds as follows. Given an
arbitrary element of $\mathrm{SD}_p$, we construct functions
$\phi$ and $\psi$ via relations~\eqref{eq:3.PhiPsiprod} and then
determine the constants~$\al_n$ by~\eqref{eq:3.al}. Since the
$\la_n$'s and $\mu_n$'s interlace, it is easily seen that all
$\al_n$'s are positive. By virtue of the results of~\cite{HMzero}
there exist functions $f_1$ and $f_2$ in $L_p(-1,1)$ such that
\begin{equation}\label{eq:5.intrepr}
\begin{aligned}
    \phi(\la) &= \cos\la + \int_{-1}^1 f_1(s)\re^{i\la s}\, ds,\\
    \psi(\la) &= \sin\la + \int_{-1}^1 f_2(s)\re^{i\la s}\, ds.
\end{aligned}
\end{equation}
Therefore the proof of Theorem~\ref{thm:3.alas} remains valid, and
thus the numbers $\al_n$ satisfy the asymptotics $\al_n = 1 +
e_{n}(g)$ for some $g\in L_p(0,1)$. We shall now prove that the
series~\eqref{eq:4.Hexp} converges in $L_p(\sM_2)$, so that the
function $H$ is well defined and belongs to~$L_p(\sM_2)$.

\begin{lemma}\label{lem:5.H}
Assume that the numbers $\la_n$ and $\al_n$, $n\in\bZ$, satisfy
the asymptotics of Theorems~\ref{thm:3.evas} and \ref{thm:3.alas}
for some $p\in[1,\infty)$. Then the series~\eqref{eq:4.Hexp}
converges in the space $L_p\bigl((-1,1),\sM_2\bigr)$ (in the
Ces\`{a}ro sense if $p=1$).
\end{lemma}

\begin{proof}
Since the matrix $B$ is skew-adjoint and has eigenvalues $\pm i$,
it suffices to prove that the scalar series
\begin{equation}\label{eq:5.h}
    h(s):= \mathrm{V.p.}\sum_{n=-\infty}^\infty
        (\al_n \re^{2\la_n is} - \re^{(2n+1)\pi is})
\end{equation}
converges in $L_p(-1,1)$ (in the Ces\`{a}ro sense if $p=1$). We
shall justify the convergence on $(0,1)$; that on $(-1,0)$ will
then follow if we replace $\al_n$ with $\al_{-1-n}$ and $\la_n$
with $-\la_{-1-n}$.

By assumption, $\al_n = 1 + e_n(g)$ for some $g\in L_p$, and
classical theorems of harmonic analysis~\cite[Sec.~I.2]{Ka} show
that the series
\[
    \mathrm{V.p.}\sum_{n=-\infty}^\infty
        (\al_n -1) \re^{2n\pi is}
\]
converges in $L_p$ to $g$ in the required sense. It remains to
establish that the series
\begin{equation}\label{eq:5.ser}
    \mathrm{V.p.}\sum_{n=-\infty}^\infty
        \al_n (\re^{2\la_n is} - \re^{(2n+1)\pi is})
\end{equation}
is convergent in $L_p$. We shall treat only the case $p=1$; the
arguments remain the same for $p>1$ but with partial Ces\`{a}ro
sums replaced by the ordinary partial sums.

By the definition of summability in the sense of Ces\`{a}ro we
have to show that the sequence of partial sums
\[
    S_N(x):= \sum_{k=1}^N \sum_{n=-k}^k
        \al_n(\re^{2\la_n ix} - \re^{(2n+1)\pi ix})
\]
is a Cauchy sequence in $L_1$. Recalling that $\la_n = \pi
(n+\tfrac12) + \tilde \la_n$ with $\tilde \la_n= e_n(f)$ for some
$f\in L_1$ and developing $\re^{2\la_nix}$ into Taylor series
around $(2n+1)\pi ix$, we see that
\begin{align*}
    S_N(x) &= \sum_{k=1}^N \sum_{n=-k}^k
        \al_n \re^{(2n+1)\pi ix} \sum_{m=1}^\infty
            \frac{(2\tilde\la_n ix)^m}{m!}\\
            &= \re^{\pi ix}\sum_{m=1}^\infty
            \frac{(2ix)^m}{m!} \sum_{k=1}^N \sum_{n=-k}^k
                \al_n\tilde\la_n^m \re^{2n\pi ix}\\
            &= \re^{\pi ix}\sum_{m=1}^\infty
                \frac{(2ix)^m}{m!} \si_N(x,f_m),
\end{align*}
where
\[
    \si_N(x,f_m) := \sum_{k=1}^N \sum_{n=-k}^k
         e_n(f_m) \re^{2n\pi ix}
\]
is the partial Ces\`{a}ro sum for the function
\[
    f_m := \underbrace{f\ast \cdots \ast f}_m +
            \underbrace{f\ast \cdots \ast f}_m \ast g,
            \qquad m\ge 1
\]
(see Appendix~\ref{app:fourier}). As is well known
(cf.~\cite[Ch.~II]{Ka}), for any $f\in L_1$, the partial
Ces\`{a}ro sums $\si_N(\cdot,f)$ converge to $f$ in $L_1$ and
$\|\si_N(\cdot,f)\|_{L_1}\le\|f\|_{L_1}$. Since
\[
    \|f_m\|_{L_1}\le (1+ \|g\|_{L_1})\|f\|^m_{L_1},
\]
the Lebesgue dominated convergence theorem shows that the limit
\[
    \lim_{N\to\infty} S_N (x)= \re^{-\pi ix}
        \sum_{m=1}^\infty \frac{(2ix)^m}{m!} f_m(x)
\]
exists in $L_1$, whence the series~\eqref{eq:5.ser} converges in
$L_1$ in the Ces\`{a}ro sense. The lemma is proved.
\end{proof}

Next we show that the GLM equation~\eqref{eq:4.GLM} and the Krein
equation~\eqref{eq:4.krein} have unique solutions that belong to
the space~$G_p(\sM_2)$. To this end it suffices to show (see
Appendix~\ref{app:krein} and Remark~\ref{rem:4.GLM} for details)
that the operator $\sI+\sH$ is positive in $\bH$, where $\sH$ is a
Wiener--Hopf operator given by
\[
    \sH \bu (x) := \int_0^1 H(x-s)\bu (s)\, ds.
\]
As a preliminary, we show that certain systems of functions form
Riesz bases in~$L_2$ or ~$\bH$.

\begin{lemma}\label{lem:5.bases}
Assume that the sequence $(\la_n)$ is as in
Theorem~\ref{thm:5.inv}. Then
\begin{itemize}
 \item[(a)] the system $(\re^{i\la_nx})_{n\in\bZ}$ forms a Riesz
basis of $L_2(-1,1)$;
 \item[(b)] the system $(\bv_n)_{n\in\bZ}$,
 with $\bv_n:=\bc_0(\cdot,\la_n)$, forms a Riesz basis of
 $\bH$.
\end{itemize}
\end{lemma}

\begin{proof}
In view of relation~\eqref{eq:3.phiint} the numbers $\la_k$ are
zeros of the exponential function of sine type with indicator
diagram~$[-1,1]$ \cite{Levin}, so that item (a) follows
from~\cite[Proposition~II.4.3]{AI}.

Assume now that $\bg := (g_1,g_2)^\mathrm{t}$ is an arbitrary
element of $\bH$. We extend $g_1$ and $g_2$ to functions $\tilde
g_1$ and $\tilde g_2$ on $(-1,1)$ in the even and odd way,
respectively, and set $\tilde g:= \tilde g_1 + i \tilde g_2$. By
(a), there exists a unique sequence $(a_n)_{n\in\bZ}$ in $\ell_2$,
for which
\[
    \tilde g = \sum_{n=-\infty}^\infty a_n \re^{i\la_nx},
\]
the series being convergent in $L_2(-1,1)$. Taking the even and
odd parts of the above equality, we arrive at the relations
\[
    \tilde g_1(x) =  \sum_{n=-\infty}^\infty a_n \cos (\la_nx),
        \qquad
   i\tilde g_2(x) = i\sum_{n=-\infty}^\infty a_n \sin (\la_nx),
\]
i.e., at the representation
\[
    \bg =   \sum_{n=-\infty}^\infty a_n \bv_n
\]
in~$\bH$. We observe that there is a constant $C$ independent of
$\tilde g$ such that
\[
    C^{-1}\sum_{n=-\infty}^\infty |a_n|^2
        \le \|\tilde g\|^2_{L_2} = 2 \|\bg\|^2_{\bH}
        \le  C \sum_{n=-\infty}^\infty |a_n|^2.
\]
Hence the system $(\bv_n)_{n\in\bZ}$ is a Riesz basis of~$\bH$.
\end{proof}

\begin{lemma}\label{lem:5.Iholds}
Assume that $\{(\la_n),(\mu_n)\}$ is an arbitrary element of
$\mathrm{SD}_p$ and that the function $H$ is constructed as
explained above. Then the corresponding operator $\sI+\sH$ is
positive in $\bH$.
\end{lemma}

\begin{proof}
We notice that the matrix $B$ is skew-adjoint and has eigenvalues
$\pm i$. It is clear that the subspaces $\ker (B \pm i) \otimes
L_2(0,1)$ are invariant subspaces of $\sH$. Therefore $\sH$ is
unitarily equivalent to the direct sum $\sH_+\oplus\sH_-$, where
\[
    (\sH_\pm g)(s) = \int_0^1 h(\pm(x-s))g(s)\,ds,
    \qquad g\in L_2,
\]
and $h$ is the function of~\eqref{eq:5.h}. Thus positivity of $\sI
+ \sH$ is equivalent to that of both $I+\sH_+$ and $I+\sH_-$. We
shall prove only that the operator~$I+\sH_+$ is positive in $L_2$;
the positivity of the other one is established analogously.

Since the systems $\bigl(\re^{(2n+1)\pi is}\bigr)_{n\in\bZ}$ and
$\bigl(\re^{2\la_nis}\bigr)_{n\in\bZ}$ constitute respectively an
orthonormal and a Riesz basis of $L_2$ (cf.
Lemma~\ref{lem:5.bases} and~\cite[Sect.~II.4.2]{AI}), it is easy
to see that
\begin{align*}
     ((I+\sH_+) f,f)
        &= (f,f) + \lim_{k\to\infty}
            \sum_{n=-k}^k \bigl[\al_n |(f, \re^{2\la_nis})|^2 -
            |(f, \re^{(2n+1)\pi is})|^2\bigr]\\
        &= \sum_{n=-\infty}^\infty \al_n |(f, \re^{2\la_nis})|^2.
\end{align*}
Since the numbers $\al_n$ are uniformly bounded away from zero, we
conclude that there is a $C>0$ such that $((I + \sH_+) f,f)\ge
C\|f\|^2$ for all $f\in L_2$, i.e., the operator $I+\sH_+$ is
(uniformly) positive in $L_2$.
\end{proof}

According to Theorem~\ref{thm:B.fact2}, positivity of the operator
$\sI+\sH$ implies that the GLM equation~\eqref{eq:4.GLM} with $F$
given by~\eqref{eq:4.F} and the Krein equation~\eqref{eq:4.krein}
have unique solutions $K, \tilde R\in G_p(\sM_2)$. We have shown
in Section~\ref{sec:GLM} that in the smooth case the solution $K$
is the kernel of the transformation operator $\sI+\sK$ for the
pair $(\tilde \sA_Q, \tilde \sA_0)$ with $Q:=\tilde R(\cdot,0)JB$.
Based on this result, we treat here the general case by a limiting
procedure.

\begin{theorem}\label{thm:5.TO}
Assume that $\{(\la_n),(\mu_n)\}$ is an arbitrary element of
$\mathrm{SD}_p$ and that $H\in L_p(\sM_2)$ is a function
of~\eqref{eq:4.Hexp} constructed as explained above. Let also $K$
and $\tilde R$ be the solutions of the GLM
equation~\eqref{eq:4.GLM} and the Krein
equation~\eqref{eq:4.krein} respectively. Denote by $\sK$ the
integral operator with kernel $K$. Then there exist a unique $Q\in
L_p(\sM_2)$---namely, $Q= \tilde R(\cdot,0)JB$---such that the
operator $\sI+\sK$ is the transformation operator for the pair
$\tilde \sA_Q$ and~$\tilde\sA_{0}$.
\end{theorem}

\begin{proof}
We shall approximate the function $H$ in the norm of~$L_p(\sM_2)$
by a sequence $(H_l)_{l=1}^\infty$ of real-valued, smooth (say,
infinitely differentiable) $\sM_2$-valued functions so that the
following holds:
\begin{itemize}
\item[(a)] for every $l\in\bN$, the GLM equation~\eqref{eq:4.GLM}
    with $H$ replaced by $H_l$ has a unique solution $K_l$, and
    the corresponding integral operators~$\sK_l$ converge to~$\sK$
    as $l\to\infty$ in the uniform operator topology
    of~$\bH$;
\item[(b)] for every $l\in\bN$, there exists $Q_l\in L_p(\sM_2)$ of
    the form~\eqref{eq:1.BQ}
    such that $\sI+ \sK_l$ is a transformation operator
    for the pair $\tilde \sA_{Q_l}$ and
    $\tilde \sA_{0}$;
\item[(c)] the matrix-functions $Q_l$ converge to
    $Q:=\tilde R(\cdot,0)JB$ in~$L_p(\sM_2)$.
\end{itemize}
If (a)--(c) hold, then by Theorems~\ref{thm:2.KinGp} and
\ref{thm:2.KTO} the operators $\sI+\sK_l$ converge in
$\sG_p(\sM_2)$ (and hence in the uniform operator topology of
$\bH$) to an operator $\sI + \sK_Q$, which is the transformation
operator for the pair $(\tilde \sA_{Q},\tilde \sA_{0})$. Thus $\sK
= \sK_Q$ yielding the result. The uniqueness of $Q$ is obvious.

The details are as follows. Using $(\la_n)$ and $(\mu_n)$, we
construct the sequence of constants $\al_n$ and set
\[
     H_l(s):= \sum_{n=-l}^l
        \Bigl(\al_n \re^{-2\la_n sB} - \re^{-\pi (2n+1) sB}\Bigr)
\]
(for $p=1$, we replace $H_l$ by the corresponding partial
Ces\`{a}ro sum $\tfrac{1}{l} \sum_{k=1}^l H_k$); then  $H_l \to H$
in $L_p(\sM_2)$ as $l \to \infty$ by Lemma~\ref{lem:5.H}.

We observe that this choice of $H_l$ corresponds to setting
$\al_n=1$ and $\la_n=\pi (n+\tfrac12)$ for all $n$ with $|n|> l$,
so that in view of Lemma~\ref{lem:5.Iholds} the Wiener--Hopf
operators $\sH_l$ with symbol~$H_l$ satisfy the condition
$\sI+\sH_l>0$. Hence by Corollary~\ref{cor:B.2} and
Remark~\ref{rem:4.GLM} the GLM equation~\eqref{eq:4.GLM} with $H$
replaced by $H_l$ has a unique solution~$K_l$. This solution~$K_l$
belongs to $G_p(\sM_2)$, and hence the corresponding integral
operator~$\sK_l$ is bounded. Since the relation $H_l\to H$ in
$L_p(\sM_2)$ as $l \to \infty$ implies that $\sH_l \to \sH$
in~$\sG_p(\sM_2)$, we conclude that $\sK_l \to \sK$ as
$l\to\infty$ in the topology of the space~$\sG_p(\sM_2)$---and
thus in the uniform operator topology of $\bH$, see
Appendix~\ref{app:krein}. This establishes (a).

Moreover, by the well-known result for continuous potentials
\cite[Ch.~12.3--4]{LS} the operator $\sI+\sK_l$ is the
transformation operator for the pair $(\tilde \sA_{Q_l},\tilde
\sA_{0})$ with
    \(
    Q_l(x):= K_l(x,x)B-BK_l(x,x)
    \).
As was explained in Section~\ref{sec:GLM}, this yields the
relation $Q_l(\cdot)=\tilde R_l(\cdot,0)JB$, where $\tilde R_l$ is
the solution to the Krein equation~\eqref{eq:4.krein} for $H=H_l$.
Thus (b) is fulfilled.

To prove (c), we observe that, according to the results of
Appendix~\ref{app:krein}, the solution $\tilde R$ to the Krein
equation~\eqref{eq:4.krein} depends continuously in the norm of
the space~$G_p(\sM_2)$ on the function $H\in L_p(\sM_2)$.
Therefore the sequence $\tilde R_l(\cdot,0)$ converges in
$L_p(\sM_2)$ to the function $\tilde R(\cdot,0)$. The proof is
complete.
\end{proof}

The last step of the reconstruction procedure is to show that the
numbers $\la_n$ and $\mu_n$ we have started with are the very
eigenvalues of the operators $\sA_1$ and $\sA_2$ with the
potential~$Q$ just found.

Since the solution $K$ to the GLM equation~\eqref{eq:4.GLM}
generates a transformation operator $\sI+\sK$ for the pair
$(\tilde \sA_Q, \tilde \sA_0)$, the functions $\bc(\cdot,\la):=
(\sI+\sK)\bc_0(\cdot,\la)$ belong to $\dom\ti\sA_Q$ and satisfy
the relation $\ti\sA_Q\bc(\cdot,\la)=\la\bc(\cdot,\la)$. We set
$\bc_k:=\bc(\cdot,\la_k)$ and show that these functions are
orthogonal and that the $\al_k$ are the corresponding norming
constants.

\begin{lemma}\label{lem:5.orth}
The system of  functions $\{\bc_k\}_{k\in\bZ}$ is an orthogonal
basis of $\bH$. Moreover, for the above numbers $\al_k$ (defined
at the beginning of this section), we have
\[
    \langle\bc_k,\bc_l\rangle = \al_k^{-1}\delta_{kl},
\]
where $\delta_{kl}$ is the Kronecker delta.
\end{lemma}

\begin{proof}
Denoting by $\bv_n$ the function $\bc_0(\cdot,\la_n)$ and
recalling that the integral operator~$\sF$ with kernel $F$
of~\eqref{eq:4.F} is related to~$\sK$ by
$(\sI+\sK)(\sI+\sF)(\sI+\sK)^*=\sI$ (see Appendix~\ref{app:krein}
for details), we conclude that
\[
    \langle\bc_k,\bc_l\rangle = \langle(\sI+\sK)^*(\sI+\sK)\bv_k,\bv_l\rangle
        = \langle(\sI+\sF)^{-1}\bv_k,\bv_l\rangle.
\]
Reverting the arguments of Section~\ref{sec:GLM}, we see that the
operator $\sI+\sF$ is equal to
\[
     \slim_{k\to\infty}\sum_{n=-k}^k \al_n \langle\cdot, \bv_n\rangle\bv_n.
\]
Since the sequence $(\bv_n)_{n\in\bZ}$ is a Riesz basis of $\bH$
in view of Lemma~\ref{lem:5.bases}, the inverse of $\sI+\sF$ can
be represented as
\[
    (\sI+\sF)^{-1} = \slim_{n\to\infty} \sum_{m=-n}^n
            \al_m^{-1} \langle\cdot, \tilde\bv_m \rangle\tilde\bv_m,
\]
where $(\tilde\bv_m)$ is a basis biorthogonal to $(\bv_m)$
(see~\cite[Ch.~VI]{GK1}). Therefore
\[
    \langle(\sI+\sF)^{-1}\bv_k,\bv_l\rangle
        = \slim_{n\to\infty}  \sum_{m=-n}^n
            \al_m^{-1} \langle \bv_k ,\ti\bv_m\rangle
                       \langle\tilde\bv_m,\bv_l\rangle
        = \sum_{m=-\infty}^\infty\al_m^{-1} \de_{km}\de_{ml}
        = \al_k^{-1} \de_{kl}.
\]
Completeness of the system $\{\bc_n\}_{n\in\bZ}$ immediately
follows from the fact that the system $\{\bv_n\}_{n\in\bZ}$ is
complete and that $\sI+\sK$ is a homeomorphism of~$\bH$, and the
lemma is proved.
\end{proof}

Next we show that the numbers $\la_n$ are indeed the eigenvalues
of the operator $\sA_1$. According to what was said above, it
suffices to show that $c_1(1,\la_k)=0$. For the operator $\tilde
\sA_Q$, one has the Green formula
\[
    \langle\tilde \sA_Q \bc(\cdot,\la),
             \bc(\cdot,\mu)\rangle-
    \langle(\bc(\cdot,\la), \tilde \sA_Q
             \bc(\cdot,\mu)\rangle
    =c_2(1,\la)c_1(1,\mu) - c_1(1,\la)c_2(1,\mu);
\]
taking therein $\la=\la_k$ and $\mu=\la_n$ and using the previous
lemma, we arrive at the relation
\begin{equation}\label{eq:5.green}
    c_2(1,\la_k)c_1(1,\la_n) = c_1(1,\la_k)c_2(1,\la_n).
\end{equation}

Assume first that none of the numbers $c_1(1,\la_k)$ vanishes.
Relation~\eqref{eq:5.green} implies that there is a constant
$\gamma$ such that, for all $k\in \bZ$,
\[
        c_2(1,\la_k)/c_1(1,\la_k) = \gamma.
\]
Then $\bc_k$ are eigenvectors corresponding to the eigenvalues
$\la_k$ of the operator that is the restriction of $\tilde \sA_Q$
by the boundary condition $u_2(1)= \gamma u_1(1)$. In other words,
the numbers $\la_n$ are zeros of the function
\[
     c_2(1,\la) - \gamma c_1(1,\la) =
        \sin\la - \gamma\cos\la   + \int_{-1}^1 \tilde f(s)
            \re^{i\la s}\,ds
\]
for some $\tilde f\in L_p(-1,1)$. However, the standard arguments
based on Rouch\'{e}'s theorem (see, e.g., \cite[Ch.~1.3]{Ma}) show
that the zeros $\tilde\la_n$ of the function
    $c_2(1,\la) - \gamma c_1(1,\la)$ obey the different asymptotics
$\tilde\la_n= \pi n + \arctan\gamma+ \mathrm{o}(1)$, which leads
to a contradiction.

Therefore there is a $k\in\bZ$ such that $c_1(1,\la_k)=0$. Then
$c_2(1,\la_k)\ne0$ due to the uniqueness of solutions to the
equation~$\ell_Q(\bu)=\la_k\bu$, and relation~\eqref{eq:5.green}
implies that $c_2(1,\la_n)=0$ for all $n\in\bZ$. In other words,
the functions~$\bc_n$ are the eigenvectors of the operator~$\sA_1$
corresponding to the eigenvalues $\la_n$. Since by
Lemma~\ref{lem:5.orth} the system $\{\bc_n\}_{n\in\bZ}$ is
complete in $\bH$, the operator $\sA_1$ has no other eigenvalues.

It remains to prove that $\mu_n$ are the eigenvalues of $\sA_2$.
We denote by $\tilde \mu_n$ these eigenvalues and construct the
corresponding function $\tilde \psi$ of~\eqref{eq:3.PhiPsiprod}.
Recalling expression~\eqref{eq:3.al}, we conclude that
$\psi(\la_k)=\tilde \psi(\la_k)$ for all $k\in \bZ$. Since the
function $\tilde \psi$ has the representation
\[
    \tilde\psi(\la) = \sin\la +
        \int_{-1}^1 \tilde f_2(s)\re^{i\la s}\, ds
\]
for some $\tilde f_2\in L_p(-1,1)$ (see the proof of
Theorem~\ref{thm:3.evas}) and $\psi$ has a similar representation
with some $f_2\in L_p(-1,1)$ instead of $\tilde f$
by~\eqref{eq:5.intrepr}, we see that the function
$\omega:=f_2-\tilde f_2$ is such that
\[
    \int_{-1}^1 \omega(s)\re^{i\la_n s}\, ds = 0
\]
for all $n\in\bZ$. Recalling that the system of functions
$\{\re^{i\la_ns}\}_{n\in\bZ}$ is closed $L_p(-1,1)$ (this follows
from~\cite[Theorem~III]{Levins} for $p>1$ and
from~\cite[Lemma~3.3]{HMzero} for $p=1$), we conclude that $\omega
= 0$. Thus the numbers $\mu_n =\tilde \mu_n$ are the eigenvalues
of the operator~$\sA_2$, and the reconstruction procedure is
complete. As the spectral data determine the function~$H$ (and
thus the transformation operator $\sI+\sK$) unambiguously, the
potential~$Q$ is unique. The proof of Theorem~\ref{thm:5.inv} is
complete.

\medskip

We observe that, in passing, we have solved the inverse spectral
problem of reconstructing the potential~$Q$ of the Dirac operator
from the spectrum of~$\sA_1$ and the corresponding sequence of
norming constants. Namely, the following is true.

\begin{theorem}\label{thm:5.norm}
Sequences of real numbers $(\la_n)_{n\in\bZ}$ and positive numbers
$(\al_n)_{n\in\bZ}$ are respectively the sequences of eigenvalues
and norming constants of an operator $\sA_1$ for some $Q\in
L_p(\sM_2)$, $p\in[1,\infty)$, if and only if the following holds:
\begin{itemize}
\item[(i)] the numbers $\la_n$ strictly increase and obey the
    asymptotics of Theorem~\ref{thm:3.evas};
\item[(ii)] the numbers $\al_n>0$ obey the asymptotics of
Theorem~\ref{thm:3.alas}.
\end{itemize}
If (i) and (ii) hold, then $Q$ is given by $R(x,x)JB=\tilde
R(x,0)JB$, where $\tilde R$ is the solution of the Krein
equation~\eqref{eq:4.krein} with $H$ of~\eqref{eq:4.Hexp}.
\end{theorem}

The reconstruction algorithm proceeds as the previous one, except
that the first and the last step (related to the
spectrum~$(\mu_n)$) should be omitted.

In a similar manner we can also treat the inverse spectral problem
for the operator~$\sA_2$ (or operators corresponding to arbitrary
separated boundary conditions).

\medskip
\textbf{Acknowledgements.} {The authors express their gratitude to
DFG for financial support of the project \mbox{436 UKR 113/79}.
The second author gratefully acknowledges the financial support of
the Alexander von Humboldt Foundation. The second and the third
authors thank the Institute for Applied Mathematics of Bonn
University for the warm hospitality.}


\appendix

\section{Fourier transform in $L_p(0,1)$}\label{app:fourier}

For any $f\in L_p(0,1)$, we denote by $e_n(f)$, $n\in\bZ$, its
$n$-th Fourier coefficients, i.e.,
\[
    e_n(f) = \int_0^1 f(x)\re^{-2\pi nix}\, dx.
\]
We also denote by $\mathbf{e}(f)$ the sequence
$\bigl(e_n(f)\bigr){}_{n\in\bZ}$ and put
\[
    X_p := \{ \mathbf{e}(f) \mid f\in L_p(0,1) \}.
\]
The vector space $X_p$ is algebraically embedded into
$\ell_\infty(\bZ)$ and becomes a Banach space under the induced
norm
\[
    \|\mathbf{e}(f)\|_{X_p} := \|f\|_{L_p}.
\]
For any $\bx=(x_n)$ and $\by=(y_n)$ in $\ell_\infty(\bZ)$ we shall
denote by $\bx\by$ the entrywise product of $\bx$ and $\by$, i.e.,
the element of $\ell_\infty(\bZ)$ with the $n$-th entry~$x_ny_n$.

\begin{proposition}\label{pro:A.mult}
$X_p$ is a commutative Banach algebra under the entrywise
multiplication, i.e.,
\begin{equation}\label{eq:A.alg}
    \|\bx\by\|_{X_p} \le \|\bx\|_{X_p}\|\by\|_{X_p}.
\end{equation}
\end{proposition}

Indeed, inequality~\eqref{eq:A.alg} follows from the fact that
$e_n(f)e_n(g) = e_n(f\ast g)$, where
\[
    (f\ast g)(x) := \int_0^1 f(x-t) g(t)\, dt
\]
is the convolution of $f$ and $g$ ($f$ being periodically extended
to~$(-1,1)$ by $f(x+1)=f(x)$, $x\in(0,1)$), and from the
inequality
\[
    \|f\ast g\|_{L_p} \le \|f\|_{L_p}\|g\|_{L_p}.
\]

The following statement is an analogue of the well-known Wiener
lemma.

\begin{proposition}\label{pro:A.wiener}
Assume that $f\in L_p(0,1)$, where $p\in[1,\infty)$. If $1+
e_n(f)\ne0$ for all $n\in\bN$, then there exists a function $g\in
L_p(0,1)$ such that
\[
    \bigl(1 + e_n(f)\bigr)^{-1} = 1 + e_n(g).
\]
\end{proposition}

\begin{proof}
To begin with, we adjoin to $X_p$ the unit element $\de$ (with all
components equal to $1$) and denote the resulting unital algebra
by $\wh X_p$. Assume that the assumptions of the lemma hold and
denote by $x$ an element of $\wh X_p$ with components $x_n:=
1+e_n(f)$. We shall prove below that $x$ is invertible in $\wh
X_p$; since $e_n(f)\to 0$ as $n\to\infty$, this will imply that
$x^{-1} = \de + y$ for some $y\in X_p$ as required.

As is well known~\cite{Na}, the element $x$ is invertible in the
unital Banach algebra~$\wh X_p$ if and only if $x$ does not belong
to any maximal ideal of $\wh X_p$. Proceeding by contradiction,
assume that there exists a maximal ideal $\fm$ of $\wh X_p$
containing $x$. Since $\wh X_p$ includes all finite sequences and
none of $x_n$ vanishes, $\fm$ also contains all finite sequences.
Finite sequences form a dense subset of $X_p$ because the set of
all trigonometric polynomials is dense in $L_p(0,1)$. Recalling
that maximal ideals are closed, we conclude that $X_p\subset\fm$.
Next we observe that $X_p$ is a proper subset of $\fm$ (e.g., $x$
belongs to $\fm\setminus X_p$) and that $X_p$ has codimension~$1$
in $\wh X_p$. Hence $\fm = \wh X_p$, which contradicts our
assumption that $\fm$ is a maximal ideal of $\wh X_p$. As a
result, $x$ is not contained in any maximal ideal of $\wh X_p$ and
thus is invertible in~$\wh X_p$ indeed. The lemma is proved.
\end{proof}


\section{The GLM equation and factorisation of Fredholm
operators}\label{app:krein}

In this appendix, we shall explain relationships between
solubility of the GLM equation and factorisation of related
Fredholm operators in some special algebras. We refer the reader
to the books~\cite{GGK,GK} for related concepts and basic facts.

Write $\bH:=L_2\bigl((0,1),\bC^2\bigr)$ and denote by $\sB$ (by
$\sB_\infty$) the Banach algebra of all bounded (compact)
operators in $\bH$. Denote also by $P_t$, $t\in[0,1]$, the
operator in $\bH$ of multiplication by $\chi_{[0,t]}$, the
characteristic function of the interval~$[0,t]$. Set
\begin{align*}
    \sB_\infty^+ &:= \{ B\in \sB_\infty \mid \forall\, t\in [0,1],
        \quad P_tB(I-P_t) = 0\},\\
    \sB_\infty^- &:= \{ B\in \sB_\infty \mid \forall \,t\in [0,1],
        \quad (I-P_t)BP_t = 0\};
\end{align*}
then $\sB_\infty^\pm$ are closed subspaces of $\sB_\infty$ and
$\sB_\infty^+ \cap \sB_\infty^- =\{0\}$. We also observe that the
operators in $\sB_\infty^\pm$ are Volterra ones.

Recall that $\sG_p(\sM_2)$, $p\in[1,\infty)$, stands for the
algebra in $\sB_\infty$ of integral operators over $(0,1)$ with
kernels in the space~$G_p(\sM_2)$ introduced in
Section~\ref{sec:TO}. The sets
\[
    \sG^+_p(\sM_2)=\sG_p(\sM_2) \cap \sB_\infty^+,
    \qquad
    \sG^-_p(\sM_2):=\sG_p(\sM_2) \cap \sB_\infty^-,
\]
are subalgebras of $\sG_p(\sM_2)$ consisting of operators with
lower- and upper-triangular kernels respectively and
$\sG_p(\sM_2)= \sG^+_p(\sM_2) \dotplus \sG^-_p(\sM_2)$.

We say that an operator $\sI+\sL$, $\sL\in \sG_p(\sM_2)$, admits a
factorization in $\sG_p(\sM_2)$ if
\begin{equation}\label{eq:B.fact}
    \sI + \sL = (\sI+\sK^+)^{-1}(\sI+\sK^-)^{-1}
\end{equation}
with some $\sK^\pm \in \sG^\pm_p(\sM_2)$.

The following two theorems were established in~\cite{My1,My2} for
the space $L_2(0,1)$, however, their generalisation to our
situation is straightforward.

\begin{theorem}\label{thm:B.fact1}
If $\sI+\sL$ admits a factorization in $\sG_p(\sM_2)$, then the
operators $\sK^\pm=\sK^\pm(\sL)$ are unique. Moreover, the set
$\Phi_p$ of operators $\sL\in\sG_p(\sM_2)$, for which the operator
$\sI+\sL$ is factorisable, is open in $\sG_p(\sM_2)$, and the
functions
$$ \Phi_p \ni \sL\mapsto \sK^\pm(\sL)\in \sG_p(\sM_2) $$ are
continuous.
\end{theorem}

\begin{theorem}\label{thm:B.fact2}
Assume that $\sL\in \sG_p(\sM_2)$. For the operator $\sI+\sL$ to
admit a factorisation in $\sG_p(\sM_2)$, it is necessary and
sufficient that the operators $I+P_t\sL P_t$ have a trivial kernel
in $\bH$ for each $t\in[0,1]$.
\end{theorem}

We remark that for a self-adjoint operator $\sL\in\sG_p(\sM_2)$
the requirement that the operators $I+P_t\sL P_t$ have a trivial
kernel in $\bH$ for all $t\in[0,1]$ is equivalent to positivity
of~$\sI+\sL$ in $\bH$.

Assume that $\sI+\sL$ is factorisable in $\sG_p(\sM_2)$, so
that~\eqref{eq:B.fact} holds. Applying $\sI+\sK^+$ to both sides
of this equality and using the fact that
$(\sI+\sK^-)^{-1}=\sI+\tilde \sK^-$ for some $\tilde \sK^-\in
\sB_\infty^-$, we derive the relation
\begin{equation}\label{eq:B.AGLM}
    \sK^+ + \sP^+ \sL + \sP^+(\sK^+\sL) = 0,
\end{equation}
where $\sP^+$ denotes the projection operator of $\sG_p(\sM_2)$
onto $\sG^+_p(\sM_2)$ parallel to $\sG^-_p(\sM_2)$. This relation
is an abstract analogue of the Gelfand--Levitan--Marchenko (GLM)
equation; indeed, in terms of the kernels $K^+$ and $L$ of the
operators $\sK^+$ and $\sL$ we get
\begin{equation}\label{eq:B.GLM}
    K^+(x,t) + L(x,t) + \int_0^x K^+(x,s)L(s,t)\, ds =0,
    \quad 0\le t\le x\le 1,
\end{equation}
cf.~\eqref{eq:4.GLM} and \eqref{eq:4.krein}.

We see that if an operator $\sI+\sL$ is factorisable in
$\sG_p(\sM_2)$, then the abstract GLM equation~\eqref{eq:B.AGLM}
has a solution~$\sK^+ = \sK^+(\sL)\in \sG^+_p(\sM_2)$. Conversely,
if $\sK^+\in \sG^+_p(\sM_2)$ is a solution of
equation~\eqref{eq:B.AGLM}, then $\sP^+(\sK^+ + \sL +
\sK^+\sL)=0$, i.e., $\sK:=\sK^+ + \sL + \sK^+\sL$ belongs to
$\sG^-_p(\sM_2)$, and the relation
\[
    (\sI+\sK^+)(\sI+\sL) = \sI+\sK
\]
holds. Since $\sI+\sK$ has the form $(\sI+\sK^-)^{-1}$ for
$\sK^-:=(\sI+\sK)^{-1}-\sI\in \sG^-_p(\sM_2)$, we conclude that
$\sI+\sL$ is factorisable in $\sG_p(\sM_2)$.

Summing up, we derive the following assertion on the solubility of
the GLM equation~\eqref{eq:B.GLM}.

\begin{corollary}\label{cor:B.2}
Assume that $\sL$ is a selfadjoint operator in~$\sG_p(\sM_2)$ with
kernel $L$ such that $\sI+\sL$ is positive. Then
equation~\eqref{eq:B.GLM} is uniquely soluble, and the
solution~$K^+$ belongs to $G_p(\sM_2)$ and depends continuously
therein on $L\in G_p(\sM_2)$.
\end{corollary}

\end{document}